\documentclass[12pt]{amsart}
\usepackage{geometry}                
\usepackage{graphicx}
\usepackage{amssymb}
\usepackage{epstopdf}
 \usepackage{latexsym}
\usepackage{amsmath}
\usepackage{amsfonts} 
\usepackage{amsthm}

\newtheorem{Proposition}{Proposition}

\newtheorem{Corollary}[Proposition]{Corollary}
\newtheorem{Remark}{Remark}

    \def\sqr#1#2{{\vcenter{\vbox{\hrule height .#2pt
                             \hbox{\vrule width .#2pt height#1pt \kern#1pt
                                   \vrule width .#2pt}
                             \hrule height .#2pt}}}}

     \def\CC{\mathbb{C}}

    \def\RR{\mathbb{R}}
    \def\ZZ{\mathbb{Z}}

\def\be{\begin{equation}}
\def\ee{\end{equation}}

\begin{document}

\title[Scalar-type Rodrigues' formulas]{Matrix valued polynomials generated by the scalar-type Rodrigues' formulas}
\author{Rodica D. Costin}




\maketitle

\ 

\ 

\begin{abstract}
The properties of matrix valued polynomials generated by the scalar-type Rodrigues' formulas are analyzed. A general representation of these polynomials is found in terms of products of simple differential operators. The recurrence relations, leading coefficients, completeness are established, as well as, in the commutative case, the second order equations for which these polynomials are eigenfunctions and the corresponding eigenvalues, and ladder operators.

The conjecture of Dur\'an and Gr\"unbaum that if the weights are self-adjoint and positive semidefinite then they are necessarily of scalar type is proved for $Q=x$ and $Q=x^2-1$ in dimension two, and for any dimension under genericity assumptions.

Commutative classes of quasi-orthogonal polynomials are found, which satisfy all the properties usually associated to orthogonal polynomials.

\end{abstract}

\

\section{Introduction}\label{Introduction}

Orthogonal matrix valued polynomials have been the subject of substantial investigation in recent years, due not only to their intrinsic interest, but also to their applications. A beautiful introduction in the topic and a wealth of references can be found in \cite{Grunbaum}, \cite{Grun_P_T}, 
\cite{Castro_Grun}, and in the recent overview \cite{Simon}.

The present paper investigates the properties of matrix valued polynomials generated by the scalar-type Rodrigues' formulas. These have also been considered in \cite{Grunbaum_Durana} under the assumption that the weight is self-adjoint and positive semidefinite. The present paper investigates the properties of resulting polynomials under no assumptions on the "weight"; in fact such polynomials appear in classification problems of differential equations 
\cite{RDC_classification}.

It is shown that the scalar-type Rodrigues' formulas can be reformulated as a product of simple differential operators, {Proposition}\,\ref{prod_form}, and their commutation relations yield a very simple derivation of the recurrence relations, leading coefficients, completeness, and, in the commutative case, the second order equations for which these polynomials are eigenfunctions, and the corresponding eigenvalues, \S\ref{GenProp}.

The problem of finding in which cases the "weights" are indeed weights, in the sense that they are self-adjoint and positive semidefinite, is investigated in the present paper for $Q(x)=x$ and 
$Q(x)=x^2-1$. It is found that under appropriate nonresonance conditions (satisfied by generic matrices) the conjecture of Dur\'an and Gr\"unbaum in \cite{Grunbaum_Durana} is true: self-adjoint weights reduce to scalar ones ({Proposition}\,\ref{Psw}) in the sense that for some $S$ we have $W(x)=S\Lambda(x)S^{-1}$ with $\Lambda(x)$ diagonal. In two dimensions the same is true for most of the resonant cases ({Proposition}\,\ref{Psw2d}). There exists one (resonant) case when there are self-adjoint weights not reducible to scalar ones ({Proposition}\,\ref{Psw2SA}); interestingly, all the generated polynomials belong to a commutative algebra of matrices. However, these weights are positive semidefinite only when reducible to scalar ones; this completes the proof of Dur\'an-Gr\"unbaum conjecture in dimension two.

The generated polynomials in this last case form a quasi-orthogonal family which satisfies all the usual properties associated to orthogonal polynomials, and they can be generalized to higher dimensions, \S\ref{pseudo}.

\section{Main Results}

\subsection{The product representation}

In this paper $V$ denotes a finitely-dimensional complex vector space, $\mathcal{M}$ denotes the matrices $\mathcal{L}(V,V)$, $\mathcal{M}[x]$ are the $\mathcal{M}$-valued polynomials, and let $\mathcal{M}(x)$ denote the set of $\mathcal{M}$-valued functions.

Let $Q(x)$ be a polynomial degree {\em{at most}} two:
\be\label{formQ}
Q(x)=\sigma x^2+\tau x+\delta
\ee
(real or complex-valued) and let $L_{1},L_2\in\mathcal{M}$ with  $L_1$ satisfying the (nonresonance) assumption:
\be\label{assumpD1}
L_1+k\sigma\ {\mbox{is\ invertible\ in\ }}\mathcal{M}\ {\mbox{for\ all\ }} k=1,2,\ldots
\ee

Let $W(x)$ be an $\mathcal{M}$-valued function satisfying the Pearson equation
\be\label{defW}
Q(x)\, W(x)^{-1}\, W'(x)\, =\, x\, L_1\, +\, L_2
\ee

{\bf{Definition}} Let $P_n(x)$ be the $\mathcal{M}$-valued function defined by the Rodrigues formula 
\be\label{Rodrig}
P_n(x)\, =\, W(x)^{-1}\, \frac{d^n}{dx^n}\left[\, Q(x)^n\, W(x)\, \right]
\ee

Note that $P_0=I$ and $P_1(x)=(2\sigma +L_1)x+\tau +L_2$. Proposition\,\ref{prod_form} below shows that in fact $P_n(x)$ are polynomials in $x$, $P_n\in\mathcal{M}[x]$.

\begin{Remark}
For one-dimensional  $V$, $P_n(x)$ are the classical orthogonal polynomials.
\end{Remark}

\begin{Remark} 
In the construction of the present paper $V$ can also be infinite-dimensional (and then $\mathcal{M}$ is the Banach space of bounded linear operators on $V$). For simplicity we restrict here to matrices.
\end{Remark}

Relations (\ref{defW}), (\ref{Rodrig}) are used in \cite{Grunbaum_Durana} in a different order of multiplication which is, of course, not essential, since it can be changed by transposition.

The representation (\ref{defAk}), (\ref{newPn}) is new\footnote{However, a special instance was used by the author in \cite{Gen_Jacobi}.} and is instrumental in the deduction 
of the properties usually associated with orthogonal polynomials, \S\ref{GenProp}.

{\em{Notation.}} For simplicity we denote $\frac{d}{dx}=\partial_x$.

\begin{Proposition}\label{prod_form}
Denote by $\mathcal{A}_k$ the following linear operators on $\mathcal{M}(x)$
\be\label{defAk}
\mathcal{A}_k\, =\, k\, Q'(x)\, +\, x\, L_1\, +\, L_2\, +\, Q(x)\partial_x
\ee

Then the matrices $P_n$ defined by (\ref{Rodrig}) satisfy
\be\label{newPn}
 P_n=\mathcal{A}_1\mathcal{A}_2\ldots\mathcal{A}_n\, I
 \ee
 for $n\geq 1$. ($I$ denotes the identity matrix.)
\end{Proposition}

{\em{Proof.}} 

Calculating successively  the $n$ derivatives in (\ref{Rodrig}) we obtain, using (\ref{defW}) and 
(\ref{defAk}),
$$\frac{d}{dx}\left[\, Q(x)^n\, W(x)\, \right]\, =\, Q(x)^{n-1}\, W(x)\, \mathcal{A}_n\, I$$
then 
$$\frac{d^2}{dx^2}\left[\, Q(x)^n\, W(x)\, \right]\, =\, Q(x)^{n-2}\, W(x)\, \mathcal{A}_{n-1}\mathcal{A}_n\, I$$
and so on.\qed

 \subsection{The weights $W(x)$ defined by (\ref{defW}).}\label{wei} \

 (o) {\em{The case $\rm{deg}\, Q=0$}} corresponds in the one-dimensional case to the Hermite polynomials, for which the operators (\ref{defAk}) are $\mathcal{A}_k=-2x+\partial_x$. For the general matrix setting, besides the general properties proved for all $Q$, the special features of this case will not be considered in the present paper.

 (i) {\em{The case $\rm{deg}\, Q=1$}}, say $Q(x)=x$, corresponds in the one-dimensional case to the Laguerre polynomials and $L_n^{(\alpha)}$ have the representation (\ref{defAk}), (\ref{newPn}) with $\mathcal{A}_k=-x+(k+\alpha)+x\partial_x$.
 
 In the general matrix setting equation (\ref{defW}) is in this case
 \be\label{dq1}
 W'(x)\, =\, W(x)\, \left(\, \frac{1}{x}\, L_2+L_1\right)
 \ee
If the eigenvalues of $L_2$ are nonresonant, in the sense that no two eigenvalues differ by an integer, then the general solution of (\ref{dq1}) has the form $W(x)=x^{L_2}\Phi(x)$ with $\Phi(x)$ analytic \cite{Anosov-Arnold}. The behavior at infinity is exponential if no eigenvalue of $L_1$ is zero. (See \S\ref{Appe1}, \S\ref{Appe2} for details.)

 (ii) {\em{The case $\rm{deg}\, Q=2$}} corresponds in dimension one to the Jacobi polynomials. Choosing $Q(x)=x^2-1$, the Jacobi polynomials $P_k^{(\alpha,\beta)}$ have  the representation (\ref{defAk}), (\ref{newPn}) with
 \be\label{forJacobi}
 \mathcal{A}_k=(2k+\alpha+\beta)x+(\alpha-\beta)+(x^2-1)\partial_x
 \ee

 In the general matrix setting equation (\ref{defW}) is
 \be\label{twosin}
  W'(x)\, =\, W(x)\, \left(\, \frac{1}{1-x}\, A\, +\, \frac{1}{x+1}\, B\right)\,\ \ {\mbox{with}}\ A=\frac{1}{2}\left(-L_1-L_2\right),\ B=\frac{1}{2}\left(L_1-L_2\right)
  \ee
with general solution $W(x)=(1-x)^A\Phi_+(x)=(x+1)^B\Phi_-(x)$ where $\Phi_\pm$ are analytic at $x=\pm 1$ (see \S\ref{Appe1} for details).

(iii) An exotic type of polynomials  with $Q=x^2-1$ appears naturally in a problem of classification of ordinary differential equations in singular regions and was studied in \cite{Gen_Jacobi}.

(iv) {\em{The commutative cases.}} In dimension two, if $L_1$ and $L_2$ commute,  it is easy to see that one of the following cases occurs in a suitable basis of $V$: \newline
(a) $L_1$ and $L_2$ are both diagonal; \newline
(b) one of $L_i$ is a multiple of the identity (this case can be regarded as a particular instance of (a) or (c)); \newline
(c) one  of $L_i$ is a Jordan block, and the other one, an "almost" Jordan block:
$$L_i=\left(\begin{array}{cc} \lambda & 1\\ 0 & \lambda\end{array}\right)\ ,\ \ \ L_{j}=\left(\begin{array}{cc} \mu & b\\ 0 & \mu\end{array}\right)$$
Note that in this last case, as a consequence of {Proposition}\,\ref{prod_form}, all the polynomials generated by (\ref{Rodrig}) belong to the commutative algebra of matrices
\be\label{comA}
\mathcal{C}\, =\, \left\{\, \left(\begin{array}{cc} \alpha & \beta\\ 0 & \alpha\end{array}\right);\ \ \alpha,\beta\in\CC\, \right\}
\ee

In higher dimensions the commutative cases are built up of blocks as described above.

\

\subsection {{Self-adjoint weights $W(x)$.}}

Since orthogonality is relative to an interval $J$, then self-adjointness $W(x)=W(x)^*$ should be required on this interval (only); in what follows it will be assumed that  $x\in J=[0,+\infty)$ in the case $Q=x$, respectively $x\in J=[-1,1]$ for $Q=x^2-1$.

 \begin{Proposition}\label{Psw}
 
 Consider the matrix $W(x)$ satisfying the Pearson formula (\ref{defW}). 
 
 Assume the following nonresonance condition: (i) in the case $Q(x)=x$ the eigenvalues of $L_2$ do not differ by an integer, respectively, (ii) in the case $Q(x)=x^2-1$  the eigenvalues of $-L_1-L_2$, or those of $L_1-L_2$, or do not differ by an integer.
 
 If $W(x)=W(x)^*$ for $x\in J$ then $W(x)$ reduces to scalar weights. 
  \end{Proposition}
  
  The proof of Proposition\,\ref{Psw} in found in \S\ref{Pf-Psw}.
 
 \begin{Proposition}\label{Psw2d}
 
In the two-dimensional case, ${\rm{dim}}\ V=2$, the conclusions of Proposition\,\ref{Psw} hold under the more general assumption that $L_2$ is diagonalizable in case (i), respectively $-L_1-L_2$, or $L_1-L_2$, is diagonalizable in case (ii).
 
 \end{Proposition}
 
 The proof of Proposition\,\ref{Psw2d} in found in \S\ref{Pf-Psw2d}.

  \begin{Proposition}\label{Psw2SA}
 Assume ${\rm{dim}}\ V=2$.
 
The following are the only solutions of (\ref{defW}), self-adjoint for $x\in J$, not reducible to scalar weights.

Define the matrices $D,\, L,\, T$ by 
 \be\label{formTDL}
L= \left(\begin{array}{cc} \alpha & \beta\\ 0 & \alpha\end{array}\right),\ D=\left(\begin{array}{cc} \lambda & 1 \\ 0 &  \lambda\end{array}\right),\ T=  \left(\begin{array}{cc} 0 & c\\ c & d\end{array}\right),\ \ \ \ \ {\mbox{with}}\ \alpha,\beta,\lambda,c,d\in\RR,\, c\ne 0
\ee

Let $S$ be an invertible matrix and let $L_1=S^{-1}LS$ and $ L_2=S^{-1}DS$. 

The following are self-adjoint solutions of  (\ref{defW}):

(i) for $Q=x$:
\be\label{saw}
W(x)=S^*TS\, {\rm{e}}^{xL_1}x^{L_2}=   {\rm{e}}^{\alpha x}x^\lambda\, S^* \left(\begin{array}{ll} 0 & c \\ c & c(\beta x+\ln x)+d\end{array}\right)S
\ee

(ii) for $Q=x^2-1$:
\begin{multline} 
W(x)=S^*TS\, {(1-x)^A}(1+x)^B={(1-x)^{-\lambda-\alpha}}(1+x)^{\lambda-\alpha}\, S^* \left(\begin{array}{ll} 0 & c \\ c & cf(x)+d\end{array}\right)S\\
\label{saw2} {\mbox{where}}\ f(x)=   (-\beta-1)\ln({1-x})+(\beta-1)\ln(1+x)
\end{multline}

  \end{Proposition}
  
The proof of Proposition\,\ref{Psw2SA} is found in \S\ref{Pf-Psw2SA}.

Since all the weights (\ref{saw}), (\ref{saw2}) have eigenvalues of opposite sign we have the following

\begin{Corollary}

If ${\rm{dim}}\ V=2$ any self-adjoint and positive semidefinite solution of (\ref{defW}) with deg $Q>0$ reduces to scalar weights.

\end{Corollary}

 \subsection{General Properties}\label{GenProp}

\subsubsection{Commutation  relations}

The commutation relations (\ref{ide1})-(\ref{daad}) are key to the beautiful properties of the classical orthogonal polynomials:

\begin{Proposition}\label{ide12}
For $r\in\mathcal{M}(x)$ the operators $\mathcal{A}_k$ satisfy the identities
\be\label{ide1}
\mathcal{A}_k(x\, r)\, =\, x\, \mathcal{A}_k\, r\, +\, Q\, r
\ee
\be\label{ide2}
Q\, \mathcal{A}_k\, r\, =\, \mathcal{A}_{k-1}\left(Q\, r\right)
\ee
\be\label{daad}
\partial_x\mathcal{A}_k\, r\, =\, \mathcal{A}_{k+1}\partial_xr\, +\, (2\sigma k+L_1)r
 \ee
\end{Proposition}

These identities follow by an immediate calculation.

\subsubsection{Recurrence relation}

\begin{Proposition}\label{recurrence}

The matrix-valued polynomials $P_n(x)$  satisfy the following two-step recurrence relation
\be\label{recPk}
x\, P_n(x)\, =P_{n+1}(x)\alpha_n+P_n(x)\, \beta_n+P_{n-1}(x)\gamma_n
\ee
whith $\alpha_n,\beta_n,\gamma_n\in\mathcal{M}$ uniquely determined; in particular
\be\label{valal}
\alpha_n=\,\left[ L_1+(2n+1)\sigma\right]^{-1}\, \left[ L_1+(2n+2)\sigma\right]^{-1}\, \left[ L_1+(n+1)\sigma\right]
\ee

\end{Proposition}

{\em{Proof of Proposition\,\ref{recurrence}.}}

Let $q\in\mathcal{M}$. The identities (\ref{ide1}), (\ref{ide2}) give by iteration 
$$
\mathcal{A}_1\mathcal{A}_2\ldots\mathcal{A}_n(x\, q)=\mathcal{A}_1\mathcal{A}_2\ldots\mathcal{A}_{n-1}\left( x\,\mathcal{A}_n q\right)\, +\,  \mathcal{A}_1\mathcal{A}_2\ldots\mathcal{A}_{n-1}\, (Q\,q)$$
\be\label{ide3}
=\, \ldots\, =x\, \mathcal{A}_1\mathcal{A}_2\ldots\mathcal{A}_n q\, +\, n\, \mathcal{A}_1\mathcal{A}_2\ldots\mathcal{A}_{n-1}\, (Q\,q)
\ee
Using (\ref{newPn})  and (\ref{ide3}) relation (\ref{recPk}) follows if we have
$$\mathcal{A}_{n}x-nQ=\mathcal{A}_n\mathcal{A}_{n+1}\alpha_n+\mathcal{A}_n\beta_n+\gamma_n$$
which expanded yields an identity of quadratic polynomials in $x$, and by identifying the coefficients we obtain the following equations for $\alpha_n,\beta_n,\gamma_n$:
\be\label{e1}
L_1+(n+1)\sigma\, =\,\left[ L_1+(2n+1)\sigma\right]\, \left[ L_1+(2n+2)\sigma\right]\, \alpha_n
\ee
\begin{multline}\label{e2}
L_2+\tau\, =\, \, (2n\sigma+L_1)\, \beta_n\\+\left[ 2(2n+1)\sigma L_2+2(n+1)\tau L_1+ 2(n+1)(2n+1)\sigma\tau+L_1L_2+L_2L_1\right]\alpha_n
\end{multline}
\begin{multline}\label{e3}
-(n-1)\delta=\gamma_n\, \\+\, \left\{ (n\tau +L_2)\left[(n+1)\tau+L_2\right]+\delta\left[ 2(n+1)\sigma+L_1\right]\right\}\,\alpha_n\, +\, (n\tau+L_2)\beta_n
\end{multline}
Due to condition  (\ref{assumpD1}), the system (\ref{e1})-(\ref{e2}) has a unique solution $\alpha_n,\beta_n,\gamma_n$ . \qed

\subsubsection{A second order equation for which $P_n$ are eigenfunctions in the commutative case.}\label{eigenS}


\begin{Proposition}\label{equation}
If $L_1$ and $L_2$ commute then $P_n$ are eigenfunctions for the operator $\mathcal{A}_1\partial_x$ and
\be\label{eigen}
\mathcal{A}_1\partial_xP_n=\, n\left[(n+1)\sigma+L_1\right]\, P_n
\ee

\end{Proposition}

{\em{Proof.}} Using the representation (\ref{newPn}) and (\ref{daad}) iteratively we obtain
$$\mathcal{A}_1\partial_xP_n=\mathcal{A}_1\partial_x\mathcal{A}_1\mathcal{A}_2\ldots\mathcal{A}_nI=(2\sigma+L_1)P_n+\mathcal{A}_1\mathcal{A}_2\partial_x\mathcal{A}_3\ldots\mathcal{A}_nI$$
$$=\ldots=\left[ (2\sigma+L_1)+\ldots+ (2n\sigma+L_1)    \right]P_n$$
which is (\ref{eigen}).\qed

 \subsubsection{The dominant coefficient.} It is clear that $P_n(x)$ are polynomials of degree at most $n$, as an immediate consequence of  the representation (\ref{newPn}). The assumption (\ref{assumpD1}) ensures that the degree of $P_n(x)$ is precisely $n$:
 
\begin{Proposition}\label{dc}
The coefficient of $x^n$ in $P_n(x)$ ($n\geq 1$) is
\be\label{domcf}
C_n=\left( L_1+2n\sigma\right)\, \left[ L_1+(2n-1)\sigma\right]\, \ldots\, \left[ L_1+(n+1)\sigma\right]
\ee
which is invertible under the assumption (\ref{assumpD1}).
\end{Proposition}

{\em{Proof.}} 
The dominant term of $P_n(x)$ is found by retaining only the dominant terms in (\ref{defAk}), and using (\ref{newPn}):
$$P_n(x)=\left( 2\sigma x+xL_1+\sigma x^2\partial_x\right)\ldots \left[ (2n-2)\sigma x+xL_1+\sigma x^2\partial_x\right]\, \left( 2n\sigma x+xL_1\right)+O\left( x^{n-1}\right)$$
and using the fact that $(\lambda x L_1+\sigma x^2\partial_x)x^k=(\lambda L_1+\sigma k)x^{k+1}$
the coefficient (\ref{domcf}) is found after a short iterative calculation. 
\qed

\subsubsection{Completeness}\label{Completeness} The sequence $\left\{ P_n\right\}_{n=0,1,2,\ldots}$ is complete in the following sense:

\begin{Proposition}\label{compl}

For any $p\in\mathcal{M}[x]$ a polynomial degree $n$, there exist $q_0,\ldots ,q_n\in\mathcal{M}$ so that
\be\label{lin_comb}
p(x)\, =\, \sum_{k=0}^n\, P_k(x)q_k
\ee
and the representation (\ref{lin_comb}) is unique.
\end{Proposition}

{\em{Proof.}} The decomposition (\ref{lin_comb}) follows easily by induction on $n$, relying on the fact that the dominant coefficients (\ref{domcf}) are invertible. \qed


\subsubsection{Orthogonality}\label{ortho}

In this section the interval $J$ is chosen so that $Q(x)W(x)$ vanishes at its endpoints. 
\begin{Proposition}\label{Proportho}
The following quasi-orthogonality relations hold:
\be\label{orr}
\int_J\, {P}_j(x)^*\, W(x)\, P_k(x)\, dx\, =\, 0\ \ {\mbox{for\ }} j<k
\ee
and
\be\label{orl}
\int_J\, { {P}}_j(x)^*\, W(x)^*\, P_k(x)\, dx\, =\, 0\ \ {\mbox{for\ }} j>k
\ee
provided that the integrals exist.

\end{Proposition}

Relations (\ref{orr}) and (\ref{orl}) are an immediate consequence of the Rodrigues formula (\ref{Rodrig}) using integration by parts (which holds for matrix multiplication).

\subsubsection{Ladder equations in the commutative case}\label{ladder}

The ladder relations for orthogonal polynomials have the form
$$L_nP_n(x)\, =\, a_n\, P_{n-1}(x)$$
where $L_n$ are first order linear differential operators and $a_n$ are constants.
These type of relations for orthogonal families of polynomials are found in \cite{Ismail} for quite general weight functions. In this section it is shown that in the commutative case, if the sequence $P_n$ satisfies a scalar type Rodrigues' formula, then the representation (\ref{defAk}), (\ref{newPn}) yields easily the ladder relations.

Indeed, consider operators of the form
\be\label{defL}
L=Ax+B+CQ(x)\partial_x,\ \ \ \ {\mbox{where}}\ A,B,C\in\CC
\ee
Note that the operators (\ref{defL}) commute with $\mathcal{A}_k$ modulo a multiple of $Q(x)$:
$$L\mathcal{A}_k=\mathcal{A}_kL+\theta_kQ(x)\ ,\ \ {\mbox{where}}\ \ \theta_k=2\sigma kC+CL_1-A$$
We then have
$$L_nP_n=L_n\mathcal{A}_1\mathcal{A}_2\ldots \mathcal{A}_n=\mathcal{A}_1L_n\mathcal{A}_2\ldots \mathcal{A}_n+\theta_1Q\mathcal{A}_2\ldots \mathcal{A}_n$$
and using 
(\ref{ide2}) we obtain further
$$=\mathcal{A}_1L_n\mathcal{A}_2\ldots \mathcal{A}_n+\theta_1\mathcal{A}_1\ldots \mathcal{A}_{n-1}Q$$
and iterating
$$=\mathcal{A}_1\mathcal{A}_2L_n\mathcal{A}_3\ldots \mathcal{A}_n+\left(\theta_1+\theta_2 \right)\mathcal{A}_1\ldots \mathcal{A}_{n-1}Q$$
until
$$=\mathcal{A}_2\ldots \mathcal{A}_{n-1}\left[\, \mathcal{A}_nL_n+\left(\theta_1+\ldots +\theta_n \right)Q\, \right]$$

The numbers $A,B,C$ are now determined so that
\be\label{eqladder}
\mathcal{A}_nL_n+\left(\theta_1+\ldots +\theta_n \right)Q=1
\ee
Expanding (\ref{eqladder}) we obtain a polynomial second degree in $x$, whose coefficients must vanish:
$$\left( n\sigma+\sigma+{L_1} \right)  \left( A+n\sigma\,C\right)=0$$
$$\left( \tau+{L_2} \right) A+ \left( 2\,n\sigma+{L_1}
 \right) B+n \left( n\sigma+\sigma+{L_1} \right) \tau\,C=0$$
 $$-\gamma\, \left( n -1\right) A+ \left( n\tau+{L_2} \right) B+n
 \left( n\sigma+\sigma+{L_1} \right) \gamma\,C=1$$
 This system has a unique solution for $n\geq 1$ under the assumption (\ref{assumpD1}):
 $$B=\left[\sigma\,{L_2}-\tau\,n\sigma-\tau\,{L_1}\right]\, G^{-1},\ \ C=\frac{1}{n}\ \left[2\,n\sigma+{L_1}\right]\, G^{-1},\ \ A=-n\sigma C$$
 where
 $$G={ n \left( 4\,\sigma\,\gamma-{
\tau}^{2} \right)  \left( n\sigma+{L_1} \right) +\gamma\,{{L_1}}^{2}+\sigma\,{{L_2}}^{2}-\tau\,{L_1}\,{L_2}}$$

 \subsection{Commutative families of quasi-orthogonal polynomials}\label{pseudo}
 
Let $L_1$ and $L_2$ belong, in a suitable basis of $V$, to the commutative algebra (\ref{comA}):
\be\label{L12C}
L_{i}=\left(\begin{array}{ll} \alpha_i & \beta_i \\ 0 & \alpha_i\end{array}\right),\ \ i=1,2
\ee
 (not both diagonal). Then $P_n(x)\in\mathcal{C}$ for all $n$. This is a genuine matrix setting, since there is no basis of $V$ in which these polynomials are diagonal - in other words, the algebra $\mathcal{C}$ is not semisimple. 

The polynomials $P_n$ satisfy the general properties of orthogonal polynomials, as shown in \S\ref{GenProp}. What is missing is orthogonality in the proper sense!

There exists, however, a signed weight (\ref{saw}), or (\ref{saw2}) with respect to which the polynomials are quasi-orthogonal.

An approach serving the function of orthogonality is the following. We will refer here to the case $Q=x$, the case $Q=x^2-1$ being similar. 

Noting that $P_n$ do not depend on the choice of the solution $W(x)$ of (\ref{defW}) let us choose a solution belonging to $\mathcal{C}$: taking $T=I$ in (\ref{saw}) (note that $S=I$) we obtain 
\be
W(x)={\rm{e}}^{xL_1}x^{L_2}={\rm{e}}^{x\alpha_1}x^{\alpha_2}\left(\begin{array}{cc} 1 & \beta_1 x+\beta_2\ln x\\ 0 & 1\end{array}\right)
\ee

Assume $\alpha_1<0$ and $\alpha_2>0$. As in {Proposition}\,\ref{Proportho} we have
\be\label{orrC}
\int_0^\infty\, {P}_j(x)\, W(x)\, P_k(x)\, dx\, =\, 0\ \ {\mbox{for\ }} j\ne k
\ee

Denote
\be\label{mome}
M_k\, =\, \int_0^\infty\, x^k\, W(x)\, dx\, 
\ee
which is obviously invertible.

Integration by parts gives 
\be\label{normC}
\int_0^\infty\, {P}_k(x)\, W(x)\, P_k(x)\, dx\, =\, (-1)^k\, C_k\ M_k
\ee
where $C_k$ is given by (\ref{domcf}), and therefore (\ref{normC}) is invertible.

Then the coefficients $q_k$ in the decomposition (\ref{lin_comb}) of any $p\in\mathcal{M}[x]$ as an expansion in $\{ P_k\}_{k\geq 0}$ can be found analytically as
\be\label{lindC}
q_k=(-1)^k\, C_k^{-1}\ M_k^{-1}\ \int_0^\infty\, p(x)\, W(x)\, P_k(x)\, dx\, 
\ee
providing another justification to the use of the term "quasi-orthogonal" family.

Note that the matrices $P_n(x)$ have classical orthogonal polynomials on the diagonal.

{\em{Generalization to higher dimensions}} can be made using the structure of commutative algebras which are not semisimple. 

Note that the algebra $\mathcal{C}$ can be re-written as the set of matrices of the form $\alpha I+\beta N$ where $N$ is a nilpotent, $N^2=0$, and $\alpha,\beta\in\CC$. In higher dimensions one can consider commutative algebras consisting of elements $\sum_{k=0}^{n-1}\alpha_kN^k$ where $N$ is nilpotent with $N^n=0$, and $\alpha_k\in\CC$.

\section{Acknowledgements}

The author is grateful for the warm and illuminating discussions with  Elena Berdysheva, Doron Lubinsky and Paul Nevai, and for the helpful and illuminating correspondence with Mourad Ismail. The author is particularly grateful for the very interesting comments of Percy Deift.

\section{Appendix}\label{Appe}

\subsection{Solutions of linear systems near regular singular points - Frobenius series}\label{Appe1}

A linear differential equation in $\CC^d$  with a regular singularity placed, say, at $x=0$ has the form $\mathbf{u}'(x)=\frac{1}{x}M(x)\mathbf{u}(x)$ where the matrix $M(x)$ is holomorphic at $x=0$. 
If the eigenvalues of $M(0)$ are nonresonant, in the sense that no two eigenvalues differ by an integer, then the fundamental solution has the form $U(x)=\Phi(x)x^{M(0)}$ with $\Phi(x)$ holomorphic at $x=0$ (see \cite{Anosov-Arnold}, \S 2.4). The matrix $\Phi$ satisfies the differential equation $\Phi '=\frac{1}{x}\Phi\left(M-M(0)\right)$ which has formal power series solution $\Phi(x)=I+O(x)$, and this series converges.

In resonant cases the fundamental solution is also a convergent series, only it may also contain logarithms besides besides the (noninteger) powers (see also \cite{Wasov}). 

These convergent representations of solutions at a regular singularity are called Frobenius series.

In the case (i) of \S\ref{wei} the matrix $M(x)$ has the form $M(x)=xL_1+L_2$, and, since the only singularities of the linear differential equation are $x=0$ and $x=\infty$ the function $\Phi(x)$ is entire.

In the case (ii) we use $M(x)=A+\frac{1-x}{x+1}B$ to write the Frobenius series at $x=1$, respectively $M(x)=B+\frac{x+1}{1-x}A$ for the Frobenius series at $x=-1$. The point at infinity is also a regular singularity: substituting $x=1/\xi$ in the equation $(x^2-1)\mathbf{u}'=xL_1+L_2$ we obtain
$\frac{d\mathbf{u}}{d\xi}=\frac{1}{\xi}\, \frac{1}{\xi^2-1}(L_1+\xi L_2)\mathbf{u}$ and the general solution has the form $\mathbf{u}=\Psi(\xi)\xi^{-L_1}$ with $\Psi$ analytic at $\xi=0$ if the eigenvalues of $L_1$ are nonresonant.

\subsection{Local behavior of linear systems near irregular singular points}\label{Appe2}

For systems $\mathbf{u}'(x)=(L_1+\frac{1}{x}L_2)\mathbf{u}$ the point $x=\infty$ is an irregular singularity. 
Assuming the nonresonance condition that the eigenvalues of $L_1$ are distinct, and working in the basis where $L_1$ is diagonal, there exists a transformation of the type $\mathbf{u}=(I+\frac{1}{x}S)\mathbf{v}$ after which the equation becomes $\mathbf{v}'=(L_1+\frac{1}{x}D+\frac{1}{x^2}g(x))\mathbf{v}$ with $D$ diagonal and $g$ analytic at infinity. The fundamental system of solutions has the representation
$\mathbf{v}=\hat{\Psi}(1/x)x^D\exp(L_1 x)$ where $\hat{\Psi}$ is a formal power series in $\frac{1}{x}$ which is Borel summable in generic cases, on appropriate sectors in the complex plane
(see \cite{Wasov} for formal solutions, \cite{BBRS} for Borel summability in the linear case, and for nonlinear equations see \cite{OC_RC} and the references therein).


\subsection{Proof of Proposition\,\ref{Psw}}\label{Pf-Psw}

(i) {\em{The case $Q=x$.}} By the nonresonance assumption, in particular the eigenvalues of $L_2$ are distinct, therefore after a change of basis of $V$ the matrix $L_2$ can be assumed diagonal: 
$$L_2\equiv D\equiv {\rm{diag}}\, \{\lambda_1,\ldots,\lambda_d\},\ \ \ \ \ (\lambda_i\ne \lambda_j\ {\mbox{for}}\ i\ne j)$$
and the equation (\ref{defW}) defining $W$ is 
\be\label{qis1case}
W'=W\left(\frac{1}{x}D+L_1\right)
\ee
Using the structure of solutions at a regular singular point, see \S\ref{Appe1}, there is a solution of the form $W_0=x^D\Phi(x)$ with $\Phi$ analytic at $x=0$. From (\ref{qis1case}) we obtain that
\be\label{eqfi}
\Phi'+\frac{1}{x}(D\Phi-\Phi D)=\Phi L_1
\ee
which has an analytic solution $\Phi(x)=I+x\Phi_1+O(x^2)$ with $\Phi_1$ satisfying 
\be\label{ephi1}
\Phi_1+D\Phi_1-\Phi_1 D=L_1
\ee

Clearly any solution of (\ref{qis1case}) can be written as $W(x)=KW_0(x)$ for a suitable constant matrix $K$.

Assume that $W(x)$ is a solution which is self-adjoint for $x>0$. Then each term in its power expansion $W(x)=Kx^D(I+x\Phi_1+O(x^2))$ must be self-adjoint. 

For the first term we have $Kx^D=x^{\overline{D}}K^*$ which implies $K_{ii}x^{\lambda_i}=\overline{K_{ii}}x^{\overline{\lambda_i}}$ for all $x>0$ hence $K_{ii},\lambda_i\in\RR$. Also we must have $K_{ij}x^{\lambda_j}=\overline{K_{ji}}x^{\lambda_i}$ for all $x>0$ and all $i\not=j$ which implies  $K_{ij}=0$ for $i\not=j$, therefore $K$ is diagonal and real. 

The auto-adjointness of the next term, $Kx^D\Phi_1$ implies that $\Phi_1$ is diagonal and real by similar arguments. Finally, relation 
(\ref{ephi1}) implies that $L_1$ is diagonal and real, which completes the proof.

(ii) {\em{The case $Q=x^2-1$}}  follows in a similar way, by writing the series of $W(x)$ at $x=1$ if 
$L_1+L_2$ is nonresonant, respectively at $x=-1$ if $L_1-L_2$ is nonresonant.

\subsection{Proof of Proposition\,\ref{Psw2d}}\label{Pf-Psw2d} \ 

Note that if the matrices assumed diagonalizable have equal eigenvalues, then the proof is exactly as the proof of Proposition\,\ref{Psw}.

 In the other resonant cases a solution of (\ref{qis1case}) has the form $W_0=x^D\Phi(x)$, where now 
$\Phi(x)$ is a convergent series containing both $x$ and $\ln x$ \cite{Wasov}. The type of logarithmic terms is determined by the resonant terms. 

Consider for example the case 
$$L_2\equiv D\equiv {\rm{diag}}\, \{\lambda_1,\lambda_2\}\ \ \ {\mbox{with}}\ \lambda_2-\lambda_1=1$$

In this case $W_0=x^D\Phi(x)$ is a solution if
\be\label{resfi}
\Phi(x)=I+ \sum_{n\geq 1}x^n\Phi_n+\sum_{n\geq 1}c_nx^n\ln x \ E_{21}
\ee
where $E_{21}$ is the matrix with entry $1$ in the position $(2,1)$ and zero everywhere else. 

Substituting (\ref{resfi}) into (\ref{eqfi}) we obtain that $\Phi_1$  must satisfy
\be\label{ephi1r}
\Phi_1+D\Phi_1-\Phi_1 D+c_1E_{21}=L_1
\ee
which has a unique solution when $c_1$ equals the entry (2,1) of $L_1$. From this point on the proof follows as in \S\ref{Pf-Psw}.

The general case $\lambda_i-\lambda_j\in\ZZ\setminus 0$ is similar, with the position of $\ln x$ in (\ref{resfi}) modified according to the resonance - see \cite{Wasov} for a detailed description of Frobenius series.

\subsection{Proof of Proposition\,\ref{Psw2SA}}\label{Pf-Psw2SA}\ 

(i) Noting that $L_1$ and $L_2$ commute then clearly (\ref{saw}) are solutions of (\ref{defW}). It only remains to show that these are the only self-adjoint solutions not reducing to the scalar case.

Using Propositions\,\ref{Psw} and\,\ref{Psw2d} the weights not reducing to the scalar case can only be found when $L_2$ is not diagonalizable, therefore when $L_2=S^{-1}DS$ with $D$ as in 
(\ref{formTDL}),  for some $S$. 

The line of the proof is similar to the one of \S\ref{Pf-Psw} (only with a different outcome!). Details are provided below, for completeness.

Denote $L=SL_1S^{-1}$. Equation (\ref{defW}) becomes $\tilde{W}'=\tilde{W}\left(\frac{1}{x}D+L\right)$ where $\tilde{W}=WS$, therefore $\tilde{W}=K{\rm{e}}^D\Phi(x)$ with $\Phi$ satisfying (\ref{eqfi}), $\Phi(x)=I+x\Phi_1+O(x^2)$. We have $W=K{\rm{e}}^D\Phi S $. 

Denote $T=S^{-1}K$ (a constant matrix); the relation $W=W^*$ is equivalent to the self-adjointness of $T{\rm{e}}^D\Phi$ for $x>0$, which implies  self-adjointness of all the terms in the series, in particular of the first term $Tx^D$, and of the second term $Tx^D\Phi_1$.

We have 
$$x^D=x^\lambda\left(\begin{array}{ll} 1 & \ln x\\ 0 & 1\end{array}\right)$$
and a simple calculations shows that $Tx^D$ is self-adjoint only for $\lambda\in\RR$ and $T$ as in (\ref{formTDL}). 

Then $Tx^D\Phi_1$ is self-adjoint only for matrices of the form
$$\Phi_1= \left(\begin{array}{cc} \alpha & \beta\\ 0 & \alpha\end{array}\right),\ \ \ \alpha,\beta\in\RR$$

Since $\Phi_1$ must solve (\ref{ephi1}) (with $L_1=L$), it follows that $L=\Phi_1$, which completes the proof of Proposition\,\ref{Psw2SA} in the case $Q=x$.

(ii) The case $Q=x^2-1$ is similar.\qed

\end{document}